\documentclass[a4paper,11pt]{article}
\usepackage{amssymb}

\newcommand{\D}{\displaystyle}

\newcommand{\R}{\mathbb{R}}

\newcommand{\N}{\mathbb{N}}

\newcommand{\beq}{\begin{equation} }
\newcommand{\eqq}{\end{equation} }
\newcommand{\cuad}{{\sqcap\kern-.68em\sqcup}}

\newcommand{\norm}[1]{\|#1\|}

\newtheorem{teo}{Theorem}[section]

\newtheorem{proposition}{Proposition}[section]

\newtheorem{lemma}{Lemma}[section]

\newtheorem{remark}{Remark}[section]
\newcommand{\bremark}{\begin{remark} \em}
\newcommand{\eremark}{\end{remark} }

\def\beeq{\begin{equation}}
\def\eeq{\end{equation}}
\newcommand{\begeqaet}{\begin{eqnarray*}}
\newcommand{\eneqaet}{\end{eqnarray*}}

\begin{document}
\begin{center}{\bf  \large
Asymptotic behaviors of governing  equation of Gauged Sigma
\\[2mm]

model for Heisenberg ferromagnet  }
\bigskip
{\small

  {\sc  Huyuan Chen}
  \medskip

 Department of Mathematics, Jiangxi Normal University,\\
Nanchang, Jiangxi 330022, PR China\\
 Email: chenhuyuan@yeah.net\\[10pt]
 {\sc  Feng Zhou}
   \medskip

Center for PDEs, School of Mathematical Sciences, East China Normal University,\\
Shanghai 200241, PR China\\
 Email: fzhou@math.ecnu.edu.cn\\[10pt]

}

\begin{abstract}
In this note, we study    weak solutions of equation
\begin{equation}\label{eq 00.1}
\Delta u =\frac{4e^u}{1+e^u} -4\pi\sum^{N}_{i=1}\delta_{p_i}+4\pi\sum^{M}_{j=1}\delta_{q_j} \quad{\rm in}\;\; \mathbb{R}^2,
\end{equation}
where  $\{\delta_{p_i}\}_{i=1}^N$ (resp. $\{\delta_{q_j}\}_{j=1}^M$ ) are Dirac masses  concentrated at the points $p_i, i=1,\cdots, N$, (resp. $q_j, j=1,\cdots, M$)
and $N-M>1$.  The equation (\ref{eq 00.1}) represents a governing equation of gauged sigma model for Heisenberg ferromagnet. We show that it has a sequence of  solutions $u_\beta$  having behaviors  as $-\beta \ln |x|+ O(1)$ at infinity with a free parameter $\beta\in(2,2(N-M))$,
 and our concern in this paper is to study the asymptotic behavior of $b_{\beta}$  as $\beta$ approaching the extremal values $2$ and $2(N-M)$.

\end{abstract}
\end{center}

\vspace{1mm}

  \smallskip



\setcounter{equation}{0}
\section{Introduction}

Vortices appear in various planar condensed-matter systems  and have important applications in many fundamental areas of physics including superconductivity \cite{A,GL,JT}, particle physics \cite{HT}, optics \cite{BEC} and cosmology \cite{VS}.  The study of multiple charges vortex construction in gauged field theory
was studied by Taubes \cite{JT,T1,T2}, initiated the existence and asymptotic behaviors of static solutions of the sigma model.
 Later on, Schroers \cite{S} extended the classical $O(3)$ sigma model solved by  Belavin-Polyakov \cite{BP1} to incorporate an Abelian gauged field and
allow the existence of vortices of opposite local charges so that the vortices of negative local charges viewed as poles of a complex scalar field $u$
makes contribute to, but those positive local charges viewed as zero of $u$ do not affect, the total energy, although they give some magnetic
manifestation for their existence \cite{Y0}. In fact, these peculiar properties are all due to the absence of symmetry breaking and in order to obtain vortices of opposite magnetic alignments with an energy that takes account of both type of vortices, it suffices to impose a broken symmetry.
 After that,  Yang in \cite{Y2} established an Abelian field theory model that allows
 the coexistence of vortices and anti-vortices,  showed how vortices and anti-vortices with the coupling of gravity, namely,
 cosmic strings and anti-strings, can be constructed in the Abelian gauged field model.

After involving the magnetic field,  the sigma model for Heisenberg ferromagnet, would be transformed into the local $U(1)-$invariant action density,
$$\mathcal{L} =-\frac14 F_{\mu\nu}F^{\mu\nu}+\frac12 D_\mu \phi \overline{D^\mu \phi}  -\frac12(1-\vec{n} \cdot \phi)^2,$$
where $\vec{n}=(0,0,1)$, $\phi:S^2\to \R^3$ with $|\phi|=1$, $D_\mu$ is gauge-covariant derivatives on $\phi$, defined by
$$D_\mu \phi =\partial _\mu \phi +A_\mu(\vec{n} \times \phi), \quad \mu=0,1,2$$
and
$$F_{\mu\nu}=\partial_\mu A_\nu-\partial_\nu A_\mu.$$
Assuming the temporal gauge $A_0=0$, the total energy is derived as
\begin{eqnarray*}
E(\phi,A) &=& \frac12\int_{\R^2}\{(D_1\phi)^2+(D_2\phi)^2+ (1-\vec{n} \cdot \phi)^2+F_{12}^2\} \\
   &=& 4\pi|deg(\phi)|+\frac12 \int_{\R^2}\{(D_1\phi\pm \phi\times D_2\phi)^2+(F_{12}\mp (1-\vec{n} \cdot \phi))^2\},
\end{eqnarray*}
where deg$(\phi)$ represents the Brouwer's degree of $\phi$.

The related Bogomol'nyi equations could be stated as
$$
\arraycolsep=1pt\left\{
\begin{array}{lll}
 \displaystyle  D_1v+iD_2v=0,
\\[2mm]\phantom{----}
 \displaystyle  F_{12}=  \frac{2|v|^2}{1+|v|^2},
\end{array}\right.
$$
then, setting $u=\ln |v|^2$, it reduces into the following governing equation of the gauged sigma model
\begin{equation}\label{eq 1.1}
-\Delta u +\frac{4e^u}{1+e^u} =4\pi\sum^{N}_{i=1}\delta_{p_i}-4\pi\sum^{M}_{j=1}\delta_{q_j} \quad{\rm in}\;\;  \mathbb{R}^2,
\end{equation}
where $\delta_p$ is the Dirac mass concentrated at $p\in \R^2$.
This subject has been expanded extensively in recent years, see the works of  Chern-Yang \cite{CY}, Lin-Yang \cite{LY}, Yang \cite{Y1}
and the references therein.
 In particular, Yang \cite{Y0,Y1} obtained a sequence of
 solutions $u_\beta$ with the asymptotic behavior
 \begin{equation}\label{aymp 1}
 u_\beta(x)=-\beta\ln |x|+b_\beta+o(1)\ \ {\rm at\ infinity\ for\ } \beta\in(2,\, 2(N-M)),
 \end{equation}
 for some constant $b_\beta$ under the restriction that  $N-M>1$.
 When $N>1$, $M=0$ and there is only one magnetic monopole, it was proved in \cite{HH}  by using  ODE analysis  that   a radial extremal non-topological solution of (\ref{eq 1.1}) has the asymptotic behavior 
 $$ u_2(r)=-2\ln r-2\ln\ln r+O(1)\quad{\rm as}\quad r\to+\infty.$$

Our aim in this paper is to consider  the behavior of $b_\beta$ in (\ref{aymp 1})   as $\beta$ approaches the extremal points $2$ and $2(N-M)$, and more asymptotic behavior estimates for the solutions of (\ref{eq 1.1}).
For convenience of readers, we use some notations and follow some presentations of known results mainly from the book of Yang \cite{Y1} (see e.g. \cite{Y0}).

\medskip
We first introduce some auxiliary functions.  Let $\rho$ be a smooth monotone increasing function over $(0,+\infty)$ such that
$$\rho(t)=\arraycolsep=1pt\left\{
\begin{array}{lll}
\ln t,\quad &0<t\le 1/2
\\[2mm]\phantom{ }
0,\quad  &t\ge 1.
\end{array}
\right.$$
Let $\{p_i\}_{i=1}^N$ and $\{q_j\}_{j=1}^M$ are different points in $\R^2$. 
Set  $v_1(x)=2\sum_{i=1}^N\rho(\frac{|x-p_i|}{\varrho})$ and  $v_2(x)=2\sum_{j=1}^M\rho(\frac{|x-q_j|}{\varrho})$,
where $\varrho\in(0,1)$ such that any two balls of
 $$\{B_{\varrho}(p_i):\, i=1,\cdots N\}\cup \{ B_{\varrho}(q_j):\,  j=1,\cdots M\}$$ do not intersect.
We fix  a positive number $r_0\ge 4e$ large enough such that $B_{\varrho}(p_i), B_{\varrho}(q_j)\subset B_{r_0}(0)$ for $i=1,\cdots, N$ and $j=1,\cdots, M$.
  Let $\eta_0:[0,+\infty)\to [0,2]$ be a smooth,  non-increasing  function with compact support in $[0,1]$ such that $\int_0^1 \eta_0(r)r dr=1$, and we take also the notation
  $\eta_0(x) = \eta_0(|x|)$ for $x \in \R^2$.
Denote $v_3=\Gamma\ast \eta_0-c_0$, where $c_0= \int_0^1(-\ln r)\eta_0(r)r dr>0$, $\ast$ means the standard convolution operator and $\Gamma(x)=-\frac1{2\pi}\ln|x|$  is the fundamental solution of Laplacian in $\R^2$, i.e.
 $$-\Delta \Gamma=\delta_0\quad {\rm in}\;\; \mathcal{D'}(\R^2).$$
   Notice that  $v_3\le 0$ is a smooth function in $\R^2$ satisfying $- \Delta v_3=\eta_0\geq 0\ {\rm in}\  B_1(0)$ and
\begin{equation}\label{6-7}
  |v_3(x)+ \ln|x|+ c_0|\le 3r_0 |x|^{-1}\quad{\rm}\quad{\rm for}\;\; |x|\ge 2r_0.
\end{equation}
See Section 2 for the proof.

Denote
$$u_0(x)=-v_1(x)+v_2(x)+\beta v_3(x),$$
where $\beta$ is a positive number to be chosen later, then $u_0$ contains all singularities at the points $p_i,\, q_j$ and the infinity. We observe that a solution $u$ of (\ref{eq 1.1}) could be written as
$u=u_0+ v,$
 where the remainder term $v$ is a bounded solution of
\begin{equation}\label{eq 1.2}
\Delta v =\frac{4K_\beta e^{v_2+v}}{e^{v_1}+K_\beta e^{v_2+v}} -g_\beta \quad{\rm in}\;\; \mathbb{R}^2,
\end{equation}
with $K_\beta=e^{\beta v_3}$, $g_1=\sum_{i=1}^N 4\pi\delta_{p_i}-\Delta v_1$, $g_2=\sum_{j=1}^M 4\pi\delta_{q_j}-\Delta v_2$ and
\begin{equation}\label{gb}
 g_\beta =g_1-g_2+\beta \Delta v_3,
\end{equation}
which is a smooth function with compact support in $B_{r_0}(0)$ and verifies that
$$\int_{\R^2} g_\beta\, dx=2\pi[2(N-M)-\beta].$$

 Our main results on asymptotic
behavior of solutions states as follows.

 \begin{teo}\label{teo 0}
Let $N-M>1$, then for any $\beta\in (2,\, 2(N-M))$, problem (\ref{eq 1.2}) has a unique solution $v_\beta$ such that
\begin{equation}\label{1.1}
 v_\beta(x)= b_\beta+ O(|x|^{-\frac{\beta-2}{\beta-1}})\quad {\rm as}\quad |x|\to+\infty,
\end{equation}
where
the constant $b_\beta\in\R$   depends on $\beta$ satisfying
\begin{equation}\label{1.4}
  \lim_{\beta\to2(N-M)^-}\frac{b_\beta}{\ln (2(N-M)-\beta)}=1,
\end{equation}
\begin{equation}\label{1.5}
 1\leq \liminf_{\beta\to2^+}\frac{b_\beta}{\ln (\beta-2)}\leq  \limsup_{\beta\to2^+}\frac{b_\beta}{\ln (\beta-2)} =2.
 \end{equation}
\end{teo}
\medskip

An interesting phenomena in Theorem \ref{teo 0} is that, by (\ref{1.4}) and (\ref{1.5}), the asymptotic
behavior of $b_\beta$ as $\beta \to (2(N-M))^-$ (resp. $\beta \to 2^+$) is of order $\ln(2(N-M) - \beta)$ (resp.  $\ln (\beta -2)$).
Back to problem (\ref{eq 1.1}), a sequence of solutions  are constructed with the  behaviors as
$-\beta \ln |x|+O(1)$ at infinity with a free parameter
$\beta\in(2,2(N-M))$. {\it Normally, this type of solutions
are called as non-topological solutions,} for instance \cite{ALW,CFL,CL,PT} on non-topological solutions of Chern-Simon equation or systems. In particular, for $\beta\in(2, 4)$,  the author in \cite[Chapter 2]{Y1}, see also \cite{Y0},  gave an existence result of (\ref{eq 1.2}) and asserted  that the solution converges to a constant at infinity as mentioned by (\ref{aymp 1}). Therefore our result extend  the existence  of solutions for (\ref{eq 1.2}) with the free parameter $\beta$ in the range $(2,\, 2(N-M))$ and the uniqueness follows by comparison principle.

Our idea for the estimates (\ref{1.4}) and (\ref{1.5}) is to construct suitable super and sub solutions of (\ref{eq 1.2}), by the uniqueness to see the asymptotic behavior
from the super  and sub solutions. These super and sub solutions are  constructed by adding some suitable constants depending on $\beta$
to  the solution of  (\ref{eq 1.2}) with $\beta=N-M+1$.

\smallskip

The rest of this paper is organized as follows. In Section 2, we  show some estimates at infinity of the convolution function
$\Gamma\ast F$ for function $F$ satisfying $\int_{\R^2} Fdx=0$ and  sketch the proof of the existence
in Theorem \ref{teo 0}.
 Section 3 is devoted to the estimates for (\ref{1.4}) and (\ref{1.5}) by constructing super  and subsolutions.







\setcounter{equation}{0}
\section{Existence and Uniqueness }


In this section, we show the results on existence and uniqueness. To this end,
 we first claim that (\ref{6-7}) holds for $|x|>2r_0$. In fact, by direct computation,  we observe that
\begin{eqnarray*}
v_3(x)+\ln|x|+c_0&=&  \frac1{2\pi}\int_{B_1(0)}(-\ln|x-y|) \eta_0(y) dy+ \frac1{2\pi}\int_{B_1(0)}\ln|x| \eta_0(y) dy \\
   &=&  -\frac1{2\pi}\int_{B_1(0)} \ln(|x-y|/|x|) \eta_0(y) dy.
\end{eqnarray*}
 For $|x|>2r_0$ and $|y|\le r_0$ (since $|y|\le 1$), we have that
 $$ 1-\frac{r_0}{|x|}\le |x-y|/|x|\le 1+\frac{r_0}{|x|},$$
then
$$ |\ln(|x-y|/|x|)|\le \max\{\ln (1+\frac{r_0}{|x|}), \,-\ln(1-\frac{r_0}{|x|})\}\le  \frac{2r_0}{|x|},$$
and thus,
$$\big|v_3(x)+\ln|x|+c_0\big|\le \frac{2r_0}{|x|},$$
the claim is true.\hfill$\Box$\smallskip

Observe  that $K_\beta =e^{\beta v_3}$ is a positive smooth function verifying that
\begin{equation}\label{3.100}
e^{-2(N-M) c_1}|x|^{-\beta}\le K_\beta (x)\le e^{2(N-M) c_1}|x|^{-\beta}\quad {\rm for}\quad |x|\ge 2r_0,
\end{equation}
when $\beta$ varies from $2$ to $2(N-M)$ and $-c_1+\ln|x|\le v_3\le c_1+\ln|x|$ for $|x|\ge 2 r_0>2e$ by (\ref{6-7}).
Let
$$\mathbb{X}_\beta =\left\{ w:  \R^2 \to \R \; \Big | \; \norm{\nabla w}_{L^2(\R^2)}+\norm{w}_{L^2(\R^2,K_\beta dx)}<\infty\Big. \right\}.$$
Now we show the following result.

\begin{proposition}\label{pr 2.1}
Let $N-M>1$ and $2<\beta<2(N-M)$.
Then problem (\ref{eq 1.2}) has a unique solution $v_\beta\in \mathbb{X}_\beta$ and there exists $b_\beta\in\R$ such that
\begin{equation}\label{3.1}
 v_\beta(x)=b_\beta+O(|x|^{-\frac{\beta-2}{\beta-1}})\quad{\rm as}\quad |x|\to \infty.
\end{equation}

\end{proposition}

To prove this result, we start the analysis by doing the decay estimates at infinity.

\begin{lemma}\label{lm 2.1}
Suppose that $\Gamma$ is the fundamental solution of $-\Delta$ in $\R^2$, $F\in L^\infty(\R^2)$ has compact support, i.e.   supp$\,F \subset \overline{B_R(0)}$ for some  $R>0$,  and satisfies that
\begin{equation}\label{2.2}
\int_{\R^2} F(x)dx=0.
\end{equation}
Then we have that 
\begin{equation}\label{e 2.3}
\norm{\Gamma\ast F}_{L^\infty(\R^2)}\le   \norm{F}_{L^1(\R^2)}+R^2\ln R\norm{F}_{L^\infty(\R^2)}
\end{equation}
and
\begin{equation}\label{e 2.1}
|\Gamma\ast F(x)|\le \frac{R}{\pi |x|} \norm{F}_{L^1(\R^2)}\quad{\rm for}\quad |x|>4R.
\end{equation}

\end{lemma}
{\bf Proof.} Since supp$F\subset \overline{B_R(0)}$ and $F\in L^\infty(\R^2)$, then $F\in L^1(\R^2)$ and for $ |x|>4R$,
\begin{eqnarray*}
|\Gamma\ast F(x)|&=& \frac{1}{2\pi}\Big|\int_{B_R(0)}  \ln|x-y| F(y)dy - \int_{B_R(0)}  \ln|x| F(y)dy\Big| \\
&= & \frac{|x|^2}{2\pi}\Big|\int_{B_\frac{R}{|x|}(0)}   \ln|e_x-z|  F(|x|z)dz\Big|
\\& \le &  \frac{|x|^2}{\pi} \int_{B_\frac{R}{|x|}(0)} |z| |F(|x|z)|dz
\\& \le & \frac{R}{\pi |x|} \norm{F}_{L^1(\R^2)},
\end{eqnarray*}
where $e_x=\frac{x}{|x|}$ and we have used (\ref{2.2}) and the fact that
$$|\ln|e_x-z||\le 2|z|\le \frac{2R}{|x|}\quad {\rm for}\quad z\in B_{ {R}/{|x|}}(0)\subset B_{1/4}(0).$$
Therefor (\ref{e 2.1}) holds. Moreover, for $|x|\le 4R$, we have that
$$|\Gamma\ast F(x)|=\frac1{2\pi}\Big|\int_{B_R(0)} F(y)\ln|x-y| dy\Big|\le  \norm{F}_{L^\infty(\R^2)} R^2 \ln R,   $$
which completes the proof of (\ref{e 2.3}).\hfill$\Box$

\medskip
Replacing  the compact support assumption for $F$ by some decay at infinity, we have the following estimate.

\begin{lemma}\label{lm 2.2}
Let $\beta\in(2,2(N-M))$,  $F\in L^1(\R^2)$ verifies (\ref{2.2}) and
\begin{equation}\label{2.4}
   |F(x)|\le c_2|x|^{-\beta},\quad  \forall\, |x|\ge 1
\end{equation}
for some $c_2\geq 1$. Then we have that
\begin{equation}\label{e 2.2}
|\Gamma\ast F(x)|\le   \frac{c_2c_3}{(\beta-2)^2} |x|^{-\frac{\beta-2}{\beta-1}}\quad{\rm for\ large}\quad |x|>4e,
\end{equation}
where $c_3>0$ depends on $ \norm{F}_{L^1(\R^2)}$, but it is independent of $c_2$ and $\beta$.

\end{lemma}
{\bf Proof.}  Since $F$ satisfies (\ref{2.2}), then  for all $|x| > 4 e$, we have that
\begin{eqnarray*}
2\pi\Gamma\ast F(x)&=&   |x|^2 \int_{\R^2}   \ln|e_x-z| F(|x|z)dz +  |x|^2\ln |x|\int_{\R^2}   F(|x|z)dz
\\&=&   |x|^2 \int_{B_{R/|x|}(0)}  \ln|e_x-z|  F(|x|z)dz +|x|^2 \int_{ B_{1/2}(e_x)}  \ln|e_x-z|  F(|x|z)dz
\\&&+|x|^2 \int_{\R^2\setminus (B_{R/|x|}(0)\cup B_{1/2}(e_x))}  \ln|e_x-z|  F(|x|z)dz
\\&=: & I_1(x)+I_2(x)+I_3(x),
\end{eqnarray*}
where    $R\in (e,\frac{|x|}4)$  will be chosen latter. Here we have that $B_{R/|x|}(0)\cap B_{1/2}(e_x)=\emptyset$ for $R\leq \frac{|x|}4$. By directly computation, we have that 
\begin{eqnarray*}
|I_1(x)| &\le &   |x|^2 \int_{B_{R/|x|}(0)} |z| |F(|x|z)|dz
\\&=&   \frac{R}{|x|} \int_{B_{R}(0)} | F(y)|dy \le \frac{R}{|x|} \norm{F}_{L^1(\R^2)}.
\end{eqnarray*}
For $z\in B_{1/2}(e_x)$, we have that $|x||z|\ge \frac12|x|>e$, then $|F(|x|z)|\le c_2|x|^{-\beta}|z|^{-\beta}$ and
\begin{eqnarray*}
|I_2(x)| &\le &  c_2 |x|^{2-\beta} \int_{B_{1/2}(e_x)} (-\ln|e_x-z|) |z|^{-\beta} dz
\\&\le & c_2 2^\beta \Big(\int_{B_{1/2}(e_x)} (-\ln|e_x-z|)   dz\Big) \, |x|^{2-\beta} 
=  c_2 2^\beta   \int_{B_{1/2}(0)} (-\ln|z|)   dz   |x|^{2-\beta} \le c_4c_2 R^{2-\beta},
\end{eqnarray*}
where $c_4 =   2^{N-M}  \big(\int_{B_{1/2}(0)} (-\ln|z|)   dz\big)$.

For $z\in \R^2\setminus (B_{R/|x|}(0)\cup B_{1/2}(e_x))$, we have that $|\ln |e_x-z||\le \ln (2+|z|)$. In fact, if $|e_x-z|\geq 1$, then it follows by the fact that $|e_x-z|\leq |e_x|+|z|=1+|z|$; if $|e_x-z|\leq 1$, we have that $|e_x-z|\geq \frac12$ for  $z\in \R^2\setminus (B_{R/|x|}(0)\cup B_{1/2}(e_x))$, then $|\ln |e_x-z||\leq \ln2\leq \ln (2+|z|)$.  Together with the fact that  $|F(|x|z)|\le c_2|x|^{-\beta}|z|^{-\beta}$,
 since $|z| \ge \frac{R}{|x|} > \frac{e}{|x|}$, and thus the integration by parts gives
\begin{eqnarray*}
|I_3(x)| &\le &  c_2  |x|^{2-\beta} \int_{\R^2\setminus  B_{R/|x|}(0)  } \ln(2+|z|)\, |z|^{-\beta}dz
\\&\le &  \frac{2\pi c_2}{\beta-2} R^{2-\beta} \ln (2+\frac{R}{|x|}) + \frac{2\pi c_2}{(\beta-2)^2} R^{2-\beta}
\\&\le &   \frac{2\pi c_2}{(\beta-2)^2}\Big(2(N-M-1)\ln 3 +1\Big)R^{2-\beta}.
\end{eqnarray*}
Thus taking $R=|x|^{\frac1{\beta-1}}$ and $|x|$ sufficient large (certainly $R\in (e,\frac{|x|}4)$ is satisfied), we have that
\begin{eqnarray*}
|\Gamma\ast F(x)| &\le &   \frac{R}{2\pi |x|} \norm{F}_{L^1(\R^2)} +\frac{c_2c_4}{2\pi}  R^{2-\beta} +\frac{ c_2}{(\beta-2)^2}\Big(2(N-M-1)\ln3+1\Big) R^{2-\beta}
\\&\le &   \frac{c_5}{(\beta-2)^2} |x|^{-\frac{\beta-2}{\beta-1}},
\end{eqnarray*}
where $c_5>0$.  This ends the proof. \hfill$\Box$\medskip

Now for $\sigma\in\R$ and $s\in\N$, we  define $W^2_{s,\sigma}$ as the closure of the set of
$C^\infty$ functions over $\R^2$ with compact supports under the norm
$$\norm{\xi}^2_{W^2_{s,\sigma}}=\sum_{|\alpha|\le s} \norm{(1+|x|)^{\sigma+|\alpha|} D^\alpha\xi}_{L^2(\R^2)}^2.$$
For more details of properties of these weighted Sobolev spaces, see e.g. \cite{Can}, \cite{McO}.
Let $C_0(\R^2)$ be the set of continuous functions on $\R^2$ vanishing at infinity.

\begin{lemma}\label{lm 2.3}\cite[Lemma 2.4.5]{Y1} The following statements hold: \\
$(i)$ If $s>1$ and $\sigma>-1$, then $W^2_{s,\sigma}\subset C_0(\R^2)$.\\
$(ii)$ For $-1<\sigma <0$, the Laplace operator $\Delta: W^2_{2,\sigma}\to W^2_{0,\sigma+2}$ is one to one and the range of $\Delta$ has the characterization
$$\Delta(W^2_{2,\sigma})=\left\{ F\in W^2_{0,\sigma+2}\;\left|\;\int_{\R^2} F dx=0 \right.\right\}.$$
$(iii)$ If $\xi\in \mathbb{X}_\beta$ and $\Delta \xi=0$, then $\xi$ is a constant.

\end{lemma}

Now we are ready to prove Proposition \ref{pr 2.1}.

\noindent{\bf Proof of Proposition \ref{pr 2.1}.} The key point to study problem (\ref{eq 1.2}) is the following equation
\begin{equation}\label{eq 2.1}
\Delta w =\frac{4K_\beta e^w}{e^{v_1}+K_\beta e^w} -h \quad{\rm in}\quad  \mathbb{R}^2,
\end{equation}
where $K_\beta$  is a positive smooth function,
$h\ge0$ is a function in $C_c^\infty(\R^2)$, i.e. with compact support such that
\begin{equation}\label{2.1}
\int_{\R^2} h(x)\,dx=2\pi[2(N-M)-\beta].
\end{equation}

{\it Existence:} As it is proved in \cite[Chapter 2, Section 2.4.2]{Y1},  problem (\ref{eq 2.1}) has a  solution $w_1$, which
 is derived by considering the critical point of the energy functional
$$I(w)=\int_{\R^2}\left\{\frac12|\nabla w|^2+4\ln (e^{v_1}+K_\beta e^w) -hw \right\} dx$$
in the admissible space
$$\mathcal{A}=\left\{w\in \mathbb{X}_\beta\left|\; \int_{\R^2} \frac{4K_\beta e^w}{e^{v_1}+K_\beta e^w}dx=\int_{\R^2} h dx\right. \right\}.  $$
Moreover, $w_1$ is a classical solution of (\ref{eq 2.1}).

The subsolution of (\ref{eq 1.2}) could be constructed as
$w_-=w_1 -\Gamma\ast (h-g)-c_6$, where $c_6>0$ is a constant such that $w_1-c_6\le0$. The supersolution is given by $w^+= \tilde w_1-\tilde w_2 -v_2$,
where $\tilde w_1$ is the solution of
$$\Delta w=\frac{4K_{\beta}e^w}{1+K_{\beta}e^w}-\tilde h,$$
$\tilde w_2=\Gamma\ast (\tilde h-\tilde g)$ with $\tilde h\ge0$ being a function in $C_c^\infty(\R^2)$ such that
$$\int_{\R^2}\tilde hdx =\int_{\R^2}\tilde gdx=2\pi(2N-\beta),$$
and $\tilde g=  (\sum_{i=1}^N 4\pi\delta_{p_i}-\Delta v_1)+\beta\Delta v_3$.
Then   a solution $v$ of  (\ref{eq 1.2}) is derived by the method of super and subsolutions. Furthermore,
$$\int_{\R^2}\frac{4K_\beta e^{v}}{e^{v_1}+K_\beta e^{v}}dx=\int_{\R^2}h dx,$$
by Lemma \ref{lm 2.3}, it is known that there is a constant $b$ such that
$$v(x)\to b\quad{\rm as}\; |x|\to+\infty.$$

{\it Uniqueness. } Assume that $w_i$ with $i=1,2$ are two solutions of problem (\ref{eq 1.2}), by Lemma 2.4.5 in \cite{Y1}, verifying that
$$w_i(x)\to b_i\quad{\rm as}\; |x|\to+\infty,$$
where we may assume that $b_1\ge b_2$.
We claim  that
\begin{equation}\label{2.1..1}
 w_1\ge w_2\quad {\rm in}\; \R^2.
\end{equation}
Otherwise, it follows by $b_1\ge b_2$, there exists $x_0\in \R^2$ such that
$$w_1(x_0)-w_2(x_0)=\min_{x\in \R^2}(w_1 -w_2)(x)<0, $$
then
$$\Delta (w_1 -w_2)(x_0)\ge0,$$
which contradicts the fact that
$$\Delta (w_1 -w_2)(x_0)=\frac{4K_\beta e^{w_1(x_0)}}{e^{v_1(x_0)-v_2(x_0)}+K_\beta e^{w_1(x_0)}}-\frac{4K_\beta e^{w_2(x_0)}}{e^{v_1(x_0)-v_2(x_0)}+K_\beta e^{w_2(x_0)}}<0.$$

Thus, $w_1\ge w_2$  in $\R^2$.  If $w_1\not=w_2$, then it implies from (\ref{2.1..1}) that
$$\int_{\R^2} g_\beta dx=\int_{\R^2}\frac{4K_\beta e^{w_1}}{e^{v_1-v_2}+K_\beta e^{w_1}}dx>\int_{\R^2}\frac{4K_\beta e^{w_2}}{e^{v_1-v_2}+K_\beta e^{w_2}}dx=\int_{\R^2} g_\beta dx,$$
which is impossible. Therefore, we have that $w_1\equiv w_2$ in $\R^2$.

We conclude that for $\beta\in(2,2(N-M))$, problem (\ref{eq 1.2}) has a unique solution $v_\beta$ and $v_\beta(x)\to b_\beta$ as $|x|\to+\infty$. Then we may rewrite that
$$v_\beta=b_\beta+\Gamma\ast \Big(\frac{4K_\beta e^{v_\beta}}{e^{v_1}+K_\beta e^{v_\beta}} -g_\beta\Big),$$
where $\displaystyle  \int_{\R^2}\Big(\frac{4K_\beta e^{v_\beta}}{e^{v_1}+K_\beta  e^{v_\beta}} -g_\beta\Big) dx=0.$
Then by applying Lemma \ref{lm 2.2},  that
$$\Big|\Gamma\ast \Big(\frac{4K_\beta e^{v_\beta}}{e^{v_1-v_2}+K_\beta e^{v_\beta}} -g_\beta\Big)(x)\Big|\le c_6|x|^{-\frac{\beta-2}{\beta-1}}, \quad \forall \, |x|>1$$
for some constant $c_6$ depending on $\beta$. This completes the proof.\hfill$\Box$

\begin{remark}\label{re 2.1}
$(i)$ Since the mapping $t\mapsto \frac{4K_\beta e^t}{e^{v_1-v_2}+K_\beta e^t}$ is increasing, the method of super and sub solutions is valid to
find out the solution. By the uniqueness and constructing  a super solution $w_1$ and a sub solution $w_2$ such that $w_1\ge w_2$,
then   the unique solution of (\ref{eq 1.2}) stays between $w_1$ and $w_2$. Furthermore, we have that
$$\int_{\R^2}\frac{4K_\beta e^{w_1}}{e^{v_1-v_2}+K_\beta e^{w_1}}\, dx\ge \int_{\R^2} g_\beta dx=2\pi[2(N-M)-\beta]$$
and
$$\int_{\R^2}\frac{4K_\beta e^{w_2}}{e^{v_1-v_2}+K_\beta e^{w_2}}\, dx\le2\pi[2(N-M)-\beta].$$

$(ii)$ From the proof of uniqueness in Proposition \ref{pr 2.1}, we conclude a type of Comparison Principle:
Let $w_1,\, w_2$ be super   and subsolutions of (\ref{eq 1.2}) respectively, verifying that $b_1\ge b_2$,
where
$$w_i(x)\to b_i\quad{\rm as}\quad |x|\to+\infty,\quad i=1,2.$$
 Then
 $$w_1\ge w_2\quad{\rm in}\quad \R^2.$$

\end{remark}

\setcounter{equation}{0}
\section{ Analysis of $b_\beta$ and $ \Gamma\ast v_\beta$ }

In this section,  we refine the estimates  of $b_\beta$ and $ \Gamma\ast v_\beta$ by constructing suitable super and sub solutions of problem
 (\ref{eq 1.2}) when $\beta$ approaches to the end points of interval $(2, 2(N-M))$. To this end, let us first set
 $$ \beta_0:= N-M+1$$
 and then (\ref{eq 1.2}) has a unique solution $v_{\beta_0}$ and
 $\int_{\R^2} g_{\beta_0} dx=2\pi[2(N-M)-\beta_0]$.

  For $\beta\in(2,\, 2(N-M))$, let $w_{\beta,\tau}=v_{\beta_0}+ \tau $  with $\tau \in \R$ being a given free parameter,   and our  super and 
  subsolutions of  (\ref{eq 1.2}) will be constructed  by varying the parameter $\tau$.

 \begin{proposition}\label{pr 3.0}
Let $\beta\in(\beta_0,2(N-M))$  and $b_\beta$ be derived by Proposition \ref{pr 2.1}, then  there exists a positive constant $c_7>0$ independent of $\beta$ such that
\begin{equation}\label{a.1}
\Big| b_\beta-\ln(2(N-M)-\beta)\Big|\le c_7.
\end{equation}

\end{proposition}
{\bf Proof.}   {\it Part 1:  subsolution.}  For $\beta\in(\beta_0,2(N-M))$, denote
$$
 \tau_{\beta,1}:=\Big(\ln\Big((2(N-M)-\beta) \frac{2\pi}{d_1}\Big)\Big)_-,
$$
where $a_-=\min\{0,\,a\}$ and $ \displaystyle d_1=  \int_{\R^2} 4K_{\beta_0} e^{v_{\beta_0}-v_1} dx.$
Let
 \begin{equation}\label{tau-0}
  w_{\beta,1 } =v_{\beta_0}+ \tau_{ \beta,1} \quad {\rm and}\quad  F_{\beta,1}=\frac{4e^{\tau_{ \beta,1} }  K_\beta e^{v_{\beta_0}}}{e^{v_1-v_2}+e^{\tau_{ \beta,1}}K_\beta e^{v_{\beta_0}}} -g_{\beta},
  \end{equation}
where $g_\beta$ is given by (\ref{gb}). In particular, we take $\tau_{\beta_0,1}=0$.
Note that
\begin{eqnarray}\label{3.300}
 -\Delta w_{\beta,1} &+&\frac{4K_\beta e^{w_{\beta,1}}}{e^{v_1-v_2}+K_\beta e^{w_{\beta,1}}} -g_{\beta}\nonumber \\
 &=& \Big(\frac{4e^{\tau_{\beta,1}}  K_\beta e^{v_{\beta_0}}}{e^{v_1-v_2}+e^{\tau_{\beta,1}}K_\beta e^{v_{\beta_0}}} -g_{\beta}\Big)-\Big(\frac{4K_{\beta_0} e^{v_{\beta_0}}}{e^{v_1-v_2}+K_{\beta_0} e^{v_{\beta_0}}}-g_{\beta_0}\Big)\label{6.1}
  \\& =:&F_{\beta,1}- F_{\beta_0,1}.\nonumber
\end{eqnarray}
Since  the mapping $\beta\mapsto K_\beta$ is decreasing, $v_2 \leq 0$, then
\begin{eqnarray*}
0 <  \int_{\R^2} \frac{4 e^{\tau_{\beta,1}}  K_\beta e^{v_{\beta_0}}}{e^{v_1-v_2}+e^{\tau_{\beta,1}} K_\beta e^{v_{\beta_0}}} dx
 &\le & e^{\tau_{\beta,1}}  \int_{\R^2} 4  K_\beta e^{v_{\beta_0}} e^{-v_1}  dx
 \\ &\le &  e^{\tau_{\beta,1}}  \int_{\R^2}4 K_{\beta_0} e^{v_{\beta_0}-v_1} dx
 \\&=& 2\pi\Big(2(N-M)-\beta\Big).
\end{eqnarray*}
From the fact that  $ \displaystyle \int_{\R^2} g_{\beta} \,dx=2\pi[2(N-M)-\beta]$, we have that
\begin{equation}\label{t-1}
    \int_{\R^2}   F_{\beta,1}  \,dx \le0.
\end{equation}

Let
$$ T_{\beta,1}:=\frac1{2\pi}\int_{\R^2} F_{\beta,1}\, dx \quad {\rm and}  \quad \tilde F_{\beta,1}:=-F_{\beta,1} + T_{\beta,1}  (-\Delta v_3),$$
then we have that $ \displaystyle \int_{\R^2}\tilde F_{\beta,1}\, dx=0.$
Obviously,   $\displaystyle \int_{\R^2}F_{\beta_0,1} dx=\int_{\R^2}\tilde F_{\beta_0,1} dx=0 $.

\smallskip

\noindent
{\it Claim 1:    There exist $\nu,\, c_8>0$    such that for any $\beta\in[\beta_0,\,2(N-M))$,
\begin{equation}\label{4.1}
 \norm{\Gamma\ast\tilde F_{\beta,1}}_{L^\infty(\R^2)}\le \nu
\end{equation}
and
\begin{equation}\label{4.2}
|\Gamma\ast \tilde F_{\beta,1}(x)|\le c_8|x|^{-\frac{\beta-2}{\beta-1}} \quad{\rm for}\; |x|>4e.
\end{equation}}
{\it Proof of Claim 1:} Since the function $t\mapsto \frac{4 e^{\tau_{\beta,1}}  K_\beta e^{t+v_{\beta_0}}}{e^{v_1-v_2}+e^{\tau_{\beta,1}} K_\beta e^{t+v_{\beta_0}}}$ is increasing and $\tau_{\beta,1}\leq 0$, then   we have that
$$\frac{4 e^{\tau_{\beta,1}}  K_\beta e^{v_{\beta_0}}}{e^{v_1-v_2}+e^{\tau_{\beta,1}} K_\beta e^{v_{\beta_0}}}\le \frac{4   K_\beta e^{v_{\beta_0}}}{e^{v_1-v_2}+  K_\beta e^{v_{\beta_0}}},$$
thus,
$$- e^{2(N-M) c_1}|x|^{-\beta}\le \tilde F_{\beta,1}(x)\le -e^{2(N-M) c_1}|x|^{-\beta},\quad\forall\, |x|> 2r_0$$
and
  $$\norm{\tilde F_{\beta,1}}_{L^\infty(\R^2)}\le 4+\norm{g_1}_{L^\infty(\R^2)}+\norm{g_2}_{L^\infty(\R^2)}+(4\pi+2)(N-M)\norm{\Delta v_3}_{L^\infty(\R^2)}:=a_0.$$
From  Lemma \ref{lm 2.2},   we have that  (\ref{4.2}) holds true for $x\in \R^N\setminus B_{4e}(0)$ and then
$\Gamma\ast\tilde F_{\beta,1}$ is bounded in  $ \R^N\setminus B_{4e}(0)$. So we only have to prove that
\begin{equation}
2\pi |\Gamma\ast \tilde F_{\beta,1}(x)|=\Big|\int_{\R^2}\ln|x-y| \tilde F_{\beta,1}(y)dy \Big|\le \nu_0,\quad {\forall}\; |x|\le 4e.
\end{equation}
In fact,  for $|x|\le 4e$, we observe that
\begin{eqnarray*}
\Big|\int_{B_{r_0}(0)}\ln|x-y|\tilde F _{\beta,1}(y)\,dy\Big| &\le &\norm{\tilde F_{\beta,1}}_{L^\infty(\R^2)} \int_{B_{r_0}(0)} |\ln|x-y||\,dy
\\&\le& \pi a_0 \Big(r_0^2(|\ln r_0|+1)+4e\Big)
\end{eqnarray*}
and
\begin{eqnarray*}
\Big|\int_{B_{r_0}^c(0)}\ln|x-y| \tilde F_{\beta,1}(y)\,dy\Big| &\le&  e^{2(N-M) c_1}\int_{B_{r_0}^c(0)} \ln(4e+|y|) |y|^{-\beta}dy
\\&\le&  e^{2(N-M) c_1}\int_{B_{r_0}^c(0)} \ln(4e+|y|) |y|^{-\beta_0}dy,
\end{eqnarray*}
which imply (\ref{4.1}), where $r_0>1$.
Thus, {\it Claim 1} holds true.

\smallskip

Now we continue to construct a subsolution.
Let
$$\underline{v}_\beta  =w_{\beta,1} +\Gamma\ast (\tilde F_{\beta,1}+ F_{\beta_0,1})-2\nu_1,$$
where $\Gamma\ast (\tilde F_{\beta,1}+ F_{\beta_0,1})-2\nu_1 \le0$
for some $\nu_1$ independent of $\beta$ by  {\it Claim 1}. Then we have that
\begin{eqnarray*}
 -\Delta \underline{v}_\beta +\frac{4K_\beta e^{\underline{v}_\beta }}{e^{v_1-v_2}+K_\beta e^{\underline{v}_\beta }}-g_\beta   &\le&
-\Delta w_{\beta,1}+ \tilde F_{\beta,1}+ F_{\beta_0,1} + \frac{4K_\beta e^{w_{\beta,1}}}{e^{v_1-v_2}+K_\beta e^{w_{\beta,1} }} -g_\beta
\\&=&T_{\beta,1}  (-\Delta v_3)
\\&\le&0,
\end{eqnarray*}
where we used (\ref{t-1}).
Then $\underline{v}_\beta $ is a subsolution of (\ref{eq 1.2}) for $\beta\in [\beta_0,\,2(N-M))$.

\medskip

\noindent{\it Part 2:   supersolution.} For $\beta\in(\beta_0,\, 2(N-M))$, we denote
$$
  \tau_{\beta,2}= \left(\ln (2(N-M)-\beta)+\ln \frac{2\pi}{d_2}\right)_-,
$$
$$
  w_{\beta,2 } =v_{\beta_0}+ \tau_{ \beta,2} \quad {\rm and}\quad  F_{\beta,2}=\frac{4e^{\tau_{ \beta,2} }  K_\beta e^{v_{\beta_0}}}{e^{v_1-v_2}+e^{\tau_{ \beta,2}}K_\beta e^{v_{\beta_0}}} -g_{\beta},
$$
where
$d_2= \displaystyle \int_{\R^2} \frac{4  K_{2(N-M)} e^{v_{\beta_0}}}{e^{v_1-v_2}+ K_{2(N-M)} e^{v_{\beta_0}}}  dx.$   Then  we derive that
\begin{eqnarray*}
 -\Delta w_{\beta,2} +\frac{4K_\beta e^{w_{\beta,2}}}{e^{v_1-v_2}+K_\beta e^{w_{\beta,2}}} -g_{\beta}   =F_{\beta,2}- F_{\beta_0,2},
\end{eqnarray*}
where $F_{\beta_0,2}=F_{\beta_0,1}$ with $\tau_{\beta_0,2}=0$.
By the decreasing monotonicity of  the function $\beta\mapsto K_\beta$, we have that
\begin{eqnarray*}
  \int_{\R^2} \frac{4 e^{\tau_{\beta,2}}  K_\beta e^{v_{\beta_0}}}{e^{v_1-v_2}+e^{\tau_{\beta,2}} K_\beta e^{v_{\beta_0}}} dx
 &\ge &    \int_{\R^2} \frac{4 e^{\tau_\beta,2}  K_{2(N-M)} e^{v_{\beta_0}}}{e^{v_1-v_2}+e^{\tau_{\beta,2}} K_{2(N-M)} e^{v_{\beta_0}}}  dx
 \\ &\ge &  e^{\tau_{\beta,2}}  \int_{\R^2} \frac{4  K_{2(N-M)} e^{v_{\beta_0}}}{e^{v_1-v_2}+ K_{2(N-M)} e^{v_{\beta_0}}}  dx
 \\&=& 2\pi[2(N-M)-\beta],
\end{eqnarray*}
which, together with $\displaystyle \int_{\R^2} g_{\beta} dx=2\pi[2(N-M)-\beta]$,  implies  that
$   \D\int_{\R^2}   F_{\beta,2}  \,dx \ge0.  $
Thus we have $T_{\beta,2}:=\frac1{2\pi}\int_{\R^2}   F_{\beta,2}\, dx\ge0$.
Let   $\tilde F_{\beta,2}:=-F_{\beta,2} + T_{\beta,2}  (-\Delta v_3).$ For the choice of $\tau_{\beta,2}$,  {\it Claim 1} holds true also, that is, there exists $\nu_2>0$ such that for $\beta\in [\beta_0,\, 2(N-M))$,
$$\norm{\tilde F_{\beta,2}}_{L^\infty(\R^2)}\le \nu_2.$$

\smallskip

   Let
$$ \bar{v}_\beta =w_{\beta,2}+\nu_1+\nu_2+\Gamma\ast (\tilde F_{\beta,2}+ F_{\beta_0,2})+\nu_3,$$
where $\nu_3=\max\{0,\, \ln \frac{d_2}{d_1}\}$.
Therefore $\underline{v}_\beta\le \bar v_\beta$ in $\; {\R^2}.$ 
Since $\nu_1+\nu_2+\Gamma\ast (\tilde F_{\beta,2}+ F_{\beta_0,2})\ge0$, then we have that
\begin{eqnarray*}
 -\Delta \bar {v}_\beta &+&\frac{4K_\beta e^{\bar {v}_\beta}}{e^{v_1-v_2}+K_\beta e^{\bar {v}_\beta }}-g_\beta   \\
&\ge & -\Delta w_{\beta,2}+\tilde F_{\beta,2} + F_{\beta_0,2} + \frac{4K_\beta e^{w_{\beta,2}}}{e^{v_1-v_2}+K_\beta e^{w_{\beta,2} }} -g_\beta
\\&=& T_{\beta,2}   (-\Delta v_3)
\\&\ge & 0,
\end{eqnarray*}
where $ T_{\beta,2} \ge 0$.
Then $\bar {v}_\beta $ is a supersolution of (\ref{eq 1.2}) for $\beta\in(\beta_0,\,2(N-M))$.

From Remark \ref{re 2.1} $(i)$, we have that for any $\beta\in (\beta_0,2(N-M))$,
\begin{equation}\label{3.2}
\underline{v}_\beta\le  v_\beta\le \bar v_\beta \quad{\rm in}\quad\R^2,
\end{equation}
which implies that (\ref{a.1}) holds. We complete the proof.\hfill$\Box$\medskip

\begin{proposition}\label{pr 3.00}
 Let $\beta\in(2, \beta_0)$ and $b_\beta$ be derived by Proposition \ref{pr 2.1}, then there exists a positive constant $c_9>0$ independent of $\beta$ such that
\begin{equation}\label{b 2}
 2\ln(\beta-2) -c_9 \leq  b_\beta  \leq    \ln(\beta-2)+c_9.
\end{equation}

\end{proposition}
{\bf Proof.} {\it Part I: subsolution. }    For   $\beta\in(2, \, \beta_0)$, we denote
$$
  \tau_{\beta,3}= \Big(2\ln(\beta-2)+\ln  (\frac{2\pi (\beta_0-2)}{d_3})\Big)_-,
$$
$$
  w_{\beta,3 } =v_{\beta_0}+ \tau_{ \beta,3} \quad {\rm and}\quad  F_{\beta,3}=\frac{4e^{\tau_{ \beta,3} }  K_\beta e^{v_{\beta_0}}}{e^{v_1-v_2}+e^{\tau_{ \beta,3}}K_\beta e^{v_{\beta_0}}} -g_{\beta},
$$
where $\displaystyle d_3=4 e^{\norm{v_{\beta_0}}_{L^\infty(\R^2)}}\Big((\beta_0-2)\int_{B_{r_0}(0)}   e^{-v_1} dx+ 2\pi\Big)$. In particular, we set   $ F_{\beta_0,3}=F_{\beta_0,1}$.
Direct computation shows that
 \begin{equation}\label{6.2}
  -\Delta w_{ \beta,3}  +\frac{4K_\beta e^{w_{ \beta,3} }}{e^{v_1-v_2}+K_\beta e^{w_{ \beta,3} }} -g_{\beta}  =F_{\beta,3 }- F_{\beta_0,1}.
 \end{equation}
Observe that
\begin{eqnarray*}
0&<& \int_{\R^2} \frac{4e^{\tau_{ \beta,3}}  K_\beta e^{v_{\beta_0}}}{e^{v_1-v_2}+e^{\tau_{ \beta,3}} K_\beta e^{v_{\beta_0}}} dx
\\&\le &4e^{\tau_{\beta,3}} e^{\norm{v_{\beta_0}}_{L^\infty(\R^2)}} \int_{\R^2} K_\beta e^{-v_1} dx
 \\ &\le &  4e^{\tau_{\beta,3}} e^{\norm{v_{\beta_0}}_{L^\infty(\R^2)}}\Big( \int_{B_{r_0}(0)}   e^{-v_1} dx+ \int_{\R^2\setminus B_{r_0}(0)} |x|^{-\beta} dx \Big)
 \\&\le&  4e^{\tau_{ \beta,3}} e^{\norm{v_{\beta_0}}_{L^\infty(\R^2)}}\Big(  \int_{B_{r_0}(0)}   e^{-v_1} dx+ \frac{2\pi}{\beta-2}  \Big)
 \\&\le&\frac{4e^{\tau_{ \beta,3}}}{\beta-2} e^{\norm{v_{\beta_0}}_{L^\infty(\R^2)}} \Big((\beta_0-2)\int_{B_{r_0}(0)}   e^{-v_1} dx+ 2\pi  \Big)
 \\&\le &  2\pi (\beta_0-2),
\end{eqnarray*}
where $e^{\tau_{ \beta,3}}\leq \frac{2\pi (\beta_0-2)}{d_3}(\beta-2)^2$.
Then we have that for $\beta\in (2,\beta_0)$,
$$  \beta-2(N-M) <\frac1{2\pi} \int_{\R^2}  F_{ \beta,3} \,dx \le  (\beta_0-2)- [2(N-M)-\beta] < 0.  $$
Thus,
\begin{equation}\label{t-3}
 T_{\beta,3}:=\frac1{2\pi}\int_{\R^2} F_{ \beta,3}\, dx\in \Big( 2-2(N-M),\, 0 \Big).
\end{equation}

\smallskip

\noindent
{\it Claim 2: }  Let   $\tilde F_{\beta,3}:=-F_{\beta,3} + T_{\beta,3}  (-\Delta v_3)$, then   there exists $\nu_3>0$ such that for any
 $\beta\in(2,\,\beta_0)$,
\begin{equation}\label{4.1.1}
  |\Gamma\ast\tilde F_{\beta,3}(x) |  \leq \nu_3,\quad\forall\,x\in\R^2.
\end{equation}

The proof of {\it Claim 2}  is postponed to the end of this proof and we continue to prove Proposition \ref{pr 3.00}.   Reset
$$  \underline{v}_\beta =w_{\beta,3}-\nu_3-\nu_1+\Gamma\ast (\tilde F_{\beta,3}+ F_{\beta_0,1}).$$
Note that  $-\nu_3-\nu_1+\Gamma\ast (\tilde F_{\beta,3}+ F_{\beta_0})\le0$ and then by (\ref{t-3}),
\begin{eqnarray*}
-\Delta \underline{v}_\beta &+&\frac{4K_\beta e^{\underline{v}_\beta}}{e^{v_1-v_2}+K_\beta e^{\underline{v}_\beta }}-g_\beta   \\&\le &
-\Delta w_{\beta,3}+ \tilde F_{\beta,3}+ F_{\beta_0,1} + \frac{4K_\beta e^{w_{\beta,3}}}{e^{v_1-v_2}+K_\beta e^{w_{\beta,3} }} -g_\beta
\\&\le & T_{\beta,3}  (-\Delta v_3)\le0.
\end{eqnarray*}
Then $\underline{v}_\beta $ is a subsolution of (\ref{eq 1.2}) for $\beta\in(2, \beta_0)$. From Lemma \ref{lm 2.2}, we have that
$$\lim_{|x|\to+\infty}\Gamma\ast (\tilde F_{\beta,3}+ F_{\beta_0,3})(x)=0, $$
and then 
\begin{equation}\label{cs 1}
 \lim_{|x|\to+\infty}\underline{v}_\beta(x)=b_{\beta_0}+2\ln(\beta-2)+\ln  (2\pi (\beta_0-1)/d_3)-\nu_3-\nu_1.
\end{equation}

\medskip

{\it Part II: supersolution. } For $\beta\in(2,\beta_0))$,
let
\begin{equation}\label{tau-3}
 \tau_{\beta,4}= \left(\ln (\beta-2)-\ln d_4\right)_-
\end{equation}
 with
$  d_4= \frac{e^{-2(N-M)c_1}}{N-M-1}\frac{   e^{b_{\beta_0}} (2r_0)^{2-\beta_0}  }{1+2(\beta_0-2)  e^{b_{\beta_0}}}.$

 Since  $\displaystyle \lim_{|x|\to+\infty}v_{\beta_0}(x)=b_{\beta_0}$,  there exists $R_0> 2
 r_0$   such that
 $$\frac12 e^{b_{\beta_0}}\le e^{v_{\beta_0}}\le 2 e^{b_{\beta_0}}, \quad \forall \; |x| \ge R_0$$
and then
 $$
 e^{\tau_{\beta,4}} K_{\beta} e^{v_{\beta_0}}\le 2(\beta_0-2) e^{b_{\beta_0}}, \quad \forall \; |x| \ge R_0,$$
for $R \ge 2R_0$,
\begin{eqnarray*}
  \int_{B_{R}(0)} \frac{4 e^{\tau_{\beta,4}}  K_\beta e^{v_{\beta_0}}}{e^{v_1-v_2}+e^{\tau_{\beta,4}} K_\beta e^{v_{\beta_0}}} dx
 &> &    \int_{B_{R}(0) \setminus B_{2r_0}(0)}  \frac{4 e^{\tau_{\beta,4}}  K_{\beta} e^{v_{\beta_0}}}{1+e^{\tau_{\beta,4}} K_{\beta} e^{v_{\beta_0}}}  dx
 \\ &\ge & \frac{ 2e^{\tau_{\beta,4}} e^{b_{\beta_0}}}{1+2(\beta_0-2)  e^{b_{\beta_0}}}  \int_{B_{R}(0) \setminus B_{2r_0}(0)} e^{-2(N-M)c_1} |x|^{-\beta}dx
 \\&\ge&\frac{e^{-2(N-M)c_1}}{d_4}\frac{ 4\pi e^{b_{\beta_0}}}{1+2(\beta_0-2)  e^{b_{\beta_0}}} (2r_0)^{2-\beta_0}\Big(1-\Big(\frac{R}{2r_0}\Big)^{2-\beta}\Big)
 \\&\ge& 4\pi(N-M-1)\Big(1-\Big(\frac{R}{2r_0}\Big)^{2-\beta}\Big),
\end{eqnarray*}
where we have used the estimate (\ref{3.100}).
Thus, passing to the limit as $R\to+\infty$, there holds 
$$ \int_{\R^2} \frac{4 e^{\tau_{\beta,4}}  K_\beta e^{v_{\beta_0}}}{e^{v_1-v_2}+e^{\tau_{\beta,4}} K_\beta e^{v_{\beta_0}}} dx\ge 4\pi(N-M-1).$$
Let $   F_{\beta,4}=\frac{4e^{\tau_{ \beta,4} }  K_\beta e^{v_{\beta_0}}}{e^{v_1-v_2}+e^{\tau_{ \beta,4}}K_\beta e^{v_{\beta_0}}} -g_{\beta},$  then   from $  \int_{\R^2} g_{\beta} dx=2\pi[2(N-M)-\beta]$, we have that
 \begin{equation}\label{t-4}
   T_{\beta,4}=\frac1{2\pi} \int_{\R^2 }   F_{\beta,4}  \,dx \ge \beta-2 >0.
 \end{equation}


Let $$ E_\beta= -F_{\beta,4}\chi_{B_{2r_0}(0)}+T_{\beta,4}(-\Delta v_3)$$
then $\displaystyle  \int_{\R^2}E_\beta\, dx=0$. From {\it Claim 1} and Lemma \ref{lm 2.1}, we obtain that
\begin{equation}\label{4.1.1-1}
 \norm{ \Gamma\ast F_{\beta_0,1}}_{L^\infty(\R^2)} \le \nu_1\quad {\rm and}\quad \norm{ \Gamma\ast E_\beta}_{L^\infty(\R^2)} \le \nu_4,
\end{equation}
where $\nu_4>0$ is independent of $\beta$.

  Denote
$$\bar v_\beta =w_{\beta}+ \nu_1+\nu_4+\nu_5+\Gamma\ast (E_\beta+F_{\beta_0,1}),$$
where $\nu_5=\ln  (2\pi (\beta_0-2)/d_3)+\ln d_4$.

Note that
$$ \nu_1+\nu_4 +\Gamma\ast (E_\beta + F_{\beta_0,1})\ge0,$$
by (\ref{t-4}), we get
\begin{eqnarray*}
-\Delta \bar  v_\beta &+&\frac{4K_\beta e^{\bar  v_\beta}}{e^{v_1-v_2}+K_\beta e^{\bar  v_\beta }} -g_\beta \\&\ge &
-\Delta w_{\beta,4}+ E_\beta+ F_{\beta_0,1} + \frac{4K_\beta e^{w_{\beta,4}}}{e^{v_1-v_2}+K_\beta e^{w_{\beta,4} }} -g_\beta
\\&=& T_{\beta,4}(-\Delta v_3)+F_{\beta,4} \chi_{\R^2\setminus B_{2r_0}(0)}
\\&\ge &0.
\end{eqnarray*}
Then  $ \bar v_\beta $ is a supersolution of (\ref{eq 1.2}) for $\beta\in(2, \beta_0)$.
By the definition of $\nu_5$, we have that
$$ \lim_{|x|\to+\infty}\underline {v}_\beta(x)\le \lim_{|x|\to+\infty}\bar {v}_\beta(x),$$
 then it infers by Remark \ref{re 2.1} $(ii)$ that $ \underline{v}_\beta\le \bar v_\beta$ in$\; \R^2.$
From Remark \ref{re 2.1} $(i)$,  for any  $\beta\in(2,\beta_1)$, we get that $\underline{v}_\beta\le v_\beta\le \bar v_\beta$ in $\; \R^2,$
which implies (\ref{b 2}).

\bigskip

{\it Finally we prove Claim 2.}   Note that  there exists $z_0>0$ such that for $\beta\in(2,\,\beta_0)$,
\begin{equation}\label{4.1.1=1}
 -z_0 \le  \Gamma\ast\tilde F_{\beta,3}(x)   \le   c_{9} (1+|x|)^{-\frac{\beta-2}{\beta+1}} +z_0,\quad\ \forall\,x\in\R^2.
\end{equation}
In fact, since $g_\beta$ has compact support, then
\begin{eqnarray*}
\norm{\tilde F_{\beta,3}}_{L^\infty(\R^2)}  &\le &     4+\norm{g_\beta}_{L^\infty(\R^2)}
\\&=&  4+\norm{g_1}_{L^\infty(\R^2)}+\norm{g_2}_{L^\infty(\R^2)} +2(N-M)\norm{\eta_0}_{L^\infty(\R^2)}
\end{eqnarray*}
and
 for $x\in \R^2\setminus B_{r_0}(0)$,
$$\frac1{c_{10}}e^{\tau_{\beta,3}}|x|^{-\beta}\le \tilde F_{\beta,3}(x)\le c_{10}e^{\tau_{\beta,3}}|x|^{-\beta},$$
where $c_{10}>1$ is independent of $\beta$ and we recall that   $e^{\tau_{\beta,3}}\leq \frac{2\pi (\beta_0-2)}{d_3}(\beta-2)^2$.\smallskip

{\it Estimates for $|x|\le 4e$. }
Observe that
$$\int_{\R^2}\ln |x-y|\, \tilde F_{\beta,3}(y)\,dy=\Big(\int_{B_{r_0}(0)}+\int_{B_{2}(x)\setminus B_{r_0}(0)}+\int_{\R^2\setminus (B_{r_0}(0)\cup B_{2}(x)) } \Big)\ln |x-y| \tilde F_{\beta,3}(y)\,dy, $$
 where
\begin{eqnarray*}
\Big|\int_{B_{r_0}(0)} \ln |x-y| \,\tilde F_{\beta,3}(y)\,dy\Big| &\le & \norm{\tilde F_{\beta,3}}_{L^\infty(\R^2)} \int_{B_{r_0}(0)}\Big| \ln |x-y| \Big|dy
\\&\le& \pi r_0^2(|\ln r_0|+1)\norm{\tilde F_{\beta,3}}_{L^\infty(\R^2)},
\end{eqnarray*}
\begin{eqnarray*}
 \Big|\int_{B_2(x)\setminus B_{r_0}(0)} \ln |x-y| \,\tilde F_{\beta,3}(y)\,dy\Big|  &\leq & \norm{\tilde F_{\beta,3}}_{L^\infty(\R^2)} \int_{B_{2}(x)} \Big| \ln |x-y| \Big|dy
 \\&\le& 4(1+\ln 2)\pi \norm{\tilde F_{\beta,3}}_{L^\infty(\R^2)}
\end{eqnarray*}
and
\begin{eqnarray*}
0\le &\displaystyle \int_{ \R^2\setminus (B_{r_0}(0)\cup B_{2}(x))} \ln |x-y| \,\tilde F_{\beta,3}(y)\,dy \\[1mm]& \displaystyle   \qquad \le    c_{10}e^{\tau_{\beta,3}} \int_{\R^2\setminus  B_{r_0}(0)}   |y|^{-\beta }\ln  |y|dy
 &\le  c_{11},
\end{eqnarray*}
where $c_{11}>1$ is independent of $\beta$.
Thus,  (\ref{4.1.1}) holds true for $|x|\le 4e$.

\smallskip

{\it Estimates for $|x|> 4e$. } This is very similar to the proof of Lemma \ref{lm 2.2}.  We rewrite
$$2\pi\Gamma\ast \tilde F_{\beta,3}(x)=: I_1(x)
+I_2(x) +I_3(x),\quad {\rm for}\; |x|> 4e.$$
Then
\begin{eqnarray*}
|I_1(x)| \le  2\frac{R}{|x|} \int_{B_{R}(0)} | \tilde F_{\beta,3}(y)|dy \le 2\pi\frac{R^3}{|x|} \norm{\tilde F_{\beta,3}}_{L^\infty(\R^2)}.
\end{eqnarray*}
For $z\in B_{1/2}(e_x)$, we have that
$0<\tilde F_{\beta,3}(|x|z)\le c_{10}e^{\tau_{\beta,3}} |x|^{-\beta}|z|^{-\beta}$ and then
\begin{eqnarray*}
|I_2(x)| &\le &  c_{10}e^{\tau_{\beta,3}} |x|^{2-\beta} \Big|\int_{B_{1/2}(e_x)} \ln|e_x-z| |z|^{-\beta} dz\Big|\le c_{12}|x|^{2-\beta},
\end{eqnarray*}
where $c_{12}>0$ is independent of $\beta$.

For $z\in \R^2\setminus (B_{R/|x|}(0)\cup B_{1/2}(e_x))$, we have that
$0<\tilde F_{\beta,3}(|x|z)\le c_{10}e^{\tau_{\beta,3}}|x|^{-\beta}|z|^{-\beta}$, then
\begin{eqnarray*}
0\le I_3(x)
 \le     2\pi c_{13} e^{\tau_{\beta,3}} \Big( \frac{ R^{2-\beta}}{\beta-2} \ln (e+\frac{R}{|x|}) + \frac{2\pi c_{13}}{(\beta-2)^2 } R^{2-\beta}\Big)
 \le   c_{14}R^{2-\beta},
\end{eqnarray*}
where $c_{13},c_{14}>0$ are independent of $\beta$.

Thus, taking $R=|x|^{\frac1{\beta+1}}$ and $|x|>4e$, we have that
\begin{eqnarray*}
2\pi\Gamma\ast \tilde F_{\beta,3}(x)   \ge   - 2\frac{R^3}{|x|} \norm{\tilde F_{\beta,3}}_{L^\infty(\R^2)} -c_{12}|x|^{2-\beta}
 \ge   -  c_{15} |x|^{-\frac{\beta-2}{\beta+1}}
\end{eqnarray*}
and
\begin{eqnarray*}
2\pi \Gamma\ast \tilde F_{\beta,3}(x)  \le  c_{15}  |x|^{-\frac{\beta-2}{\beta+1}},
 \end{eqnarray*}
where $c_{15}>0$ is independent of $\beta$.
Therefore,  (\ref{4.1.1=1})    holds true.  \hfill$\Box$

\bigskip

\noindent{\bf Proof of Theorem \ref{teo 0}.}  From Proposition \ref{pr 2.1}, problem (\ref{eq 1.2}) has a unique solution $w_\beta$, which verifies (\ref{3.1}).
Proposition \ref{pr 3.0} and Proposition \ref{pr 3.00}
show that
\begin{equation}\label{1.4+1}
  |b_\beta -\ln (2(N-M)-\beta)|\le c_7\quad{\rm for}\quad \beta\in (\beta_0,\,2(N-M))
\end{equation}
and
\begin{equation}\label{1.5+1}
 2\ln(\beta-2)- c_9 \leq b_\beta   \le  \ln(\beta-2)+ c_9\quad{\rm for}\quad \beta\in (2,\,\beta_0),
\end{equation}
which imply (\ref{1.4}) and (\ref{1.5}) respectively.\hfill$\Box$

\bigskip\bigskip

 \noindent{\small {\bf Acknowledgements:} H. Chen  is supported by NSFC (No. 11661045 and 11726614).
 F. Zhou is supported by Science and Technology Commission of Shanghai Municipality (STCSM, No. 18dz2271000) and NSFC (No.11726613 and 11431005).}

\end{document}